\documentclass[11pt]{article}
\usepackage{amssymb}
\usepackage{mathrsfs}
\usepackage{latexsym,amsmath,amssymb,amsfonts,epsfig,graphicx,cite,psfrag}
\usepackage{eepic,color,colordvi,amscd}
\usepackage{ebezier}
\usepackage{verbatim}
\usepackage{comment}

\newcommand{\qed}{\hfill $\Box $}
\newcommand{\pf}{\noindent {\bf Proof.} }

\newtheorem{theorem}{Theorem}[section]
\newtheorem{lemma}[theorem]{Lemma}

\begin{document}

\title{Rainbow perfect matchings for 4-uniform hypergraphs}

\author{Hongliang Lu\footnote{Partially supported by the National Natural
Science Foundation of China under grant No.11871391 and
Fundamental Research Funds for the Central Universities}\\
School of Mathematics and Statistics\\
Xi'an Jiaotong University\\
Xi'an, Shaanxi 710049, China\\
\medskip \\
Yan Wang\\
\medskip \\
Xingxing Yu\footnote{Partially supported by NSF grant DMS-1954134}\\
School of Mathematics\\
Georgia Institute of Technology\\
Atlanta, GA 30332, USA}

\date{}
\maketitle

\begin{abstract}

Let $n$ be a sufficiently large integer with $n\equiv 0\pmod 4$ and let $F_i \subseteq{[n]\choose 4}$ where $i\in [n/4]$. We show that if each vertex of $F_i$ is contained in more than ${n-1\choose 3}-{3n/4\choose 3}$ edges, then $\{F_1, \ldots ,F_{n/4}\}$ admits a rainbow matching, i.e., a set of $n/4$ edges consisting of one edge from each $F_i$. This generalizes a deep result of Khan on perfect matchings in 4-uniform hypergraphs.
% (when $F_1=\ldots =F_{n/4}$).

\end{abstract}

\section{Introduction}

A {\it hypergraph} is a family of subsets (called {\it edges}) of a
nonempty set whose elements are the {\it vertices} of the
hypergraph.  For a hypergraph $H$, we use $V(H)$ to denote its vertex
set, and $E(H)$ to denote its edge set, and   let $e(H):=|E(H)|$.
 We say that a hypergraph $H$ is  {\it $k$-uniform} for some
positive integer $k$ if all edges of $H$
have the same size $k$. A $k$-uniform hypergraph is also known as a
{\it $k$-graph}.

A \emph{matching} in a hypergraph $H$ is a set of pairwise disjoint
edges of $H$. %and we use $\nu(H)$ to denote the maximum size of a matching in $H$.
Finding maximum matchings in $k$-graphs
is NP-hard for $k\ge 3$, see \cite{Ka72}. Hence, it is of interest
to find good sufficient conditions for the existence of a
large matching in $k$-graphs. The most well known open problem in this
area is the following
conjecture made by Erd\H{o}s
\cite{Erdos65} in 1965:
For positive integers $k,n,t$, if $H$ is a $k$-graph of order $n$ and
 $e(H)> \max \left\{{kt-1\choose k}, {n\choose k}-{n-t+1\choose
    k}\right\}$ then $H$ has a matching of size $t$.
This bound on $e(H)$ is tight because of the complete $k$-graph on $kt-1$ vertices and the $k$-graph on $n$ vertices in which every
edge intersects a fixed set of $t-1$ vertices. There have been recent
activities on this conjecture, see \cite{AHS,AFH12,FLM,Fr13,Fr17,FK-20,HLS,LM}. %\cite{AHS, AFH12, FLM, Fr13, Fr17, FrK19+, HLS, LM}.

One type of conditions that have been used to ensure the existence of
large matchings is the so-called
``Dirac-type conditions'', which involves degrees of sets of vertices. Our work in this
paper falls into this category. For convenience, let
$[k]:=\{1,\ldots, k\}$ for any positive integer $k$, and let
${S\choose k}:=\{T\subseteq S: |T|=k\}$ for any set $S$ and positive
integer $k$. Let $H$ be a hypergraph. For any $T\subseteq V(H)$, the {\it degree}
of $T$ in $H$, denoted by $d_H(T)$, is the number of edges of $H$
containing $T$. For any integer  $l\ge 0$, $\delta_l(H):=
\min\{d_H(T): T\in {V(H)\choose l}\}$ denotes the minimum
{\it $l$-degree} of $H$. Hence,  $\delta_0(H)=e(H)$.  Note that
$\delta_1(H)$ is often called the {\it minimum
vertex degree} of $H$.

For integers $n,k,s,d$ satisfying $0 \leq d \leq k-1$, $n\equiv 0 \pmod
k$,  and $0\leq s\leq n/k$, $m^s_
d(k,n)$ denotes the minimum integer $m$ such that every $k$-graph $H$ on $n$ vertices with $\delta_d (H)\geq m$
has a matching of size $s$.
R\"{o}dl, Ruci\'{n}ski, and Szemer\'{e}di \cite{RRS06} determined
$m^{n/k}_{k-1}(k,n)$ for large $n$, which has motivated a large amount
of work, see \cite{TZ13, KOT13, Kh13, Kh16, HKP18, Han19}. For instance,  Treglown and Zhao \cite{TZ13}
extended this result by determining $m^{n/k}_{d}(k,n)$ for all
$d\geq k/2$. On the other hand, it seems more difficult to determine $m^{n/k}_{d}(k,n)$
when $d<k/2$. K\"uhn, Osthus, and Treglown \cite{KOT13} and, independently, Khan
\cite{Kh13} determined $m_1^{n/3}(3,n)$.  Khan \cite{Kh16} further
determined $m_1^{n/4}(4,n)$. The main work in this paper is to prove a
more general result which implies Khan's result and uses
different techniques.

Let $\mathcal{F} = \{F_1,\ldots, F_t\}$ be a family
of hypergraphs. A set of $t$ pairwise disjoint edges, one from
each $F_i$, is called a \emph{rainbow matching} for $\mathcal{F}$. (In this case, we also say that ${\cal F}$ or
$\{F_1,\ldots,  F_t\}$ {\it admits} a rainbow matching.) There has been a lot of interest in studying rainbow versions of
matching problems, see \cite{AH09, AH17,AF,AH,FK-20,JK,KK,KL,HLS,HZ17,LYY2, MN,PY}. For instance,  Aharoni and
 Howard \cite{AH} made the following conjecture, which first appeared
 in  Huang, Loh, and Sudakov \cite{HLS}:
%\begin{conjecture}\label{HLSAH}
Let $t$ be a positive integer and ${\cal F}=\{F_1,\ldots, F_t\}$ such that, for $i\in [t]$,
$F_i\subseteq {[n]\choose k}$ and
$e(F_i)> \left\{{kt-1\choose k}, {n\choose k}-{n-t+1\choose k}\right\};$
 then  ${\cal F}$ admits a rainbow matching.
%\end{conjecture}
Huang, Loh, and Sudakov \cite{HLS} showed that this conjecture  holds
for $n> 3k^2t$. Frankl and Kupavskii \cite{FK-20} improved this lower
bound to
$n\ge 12tk\log(e^2t)$, which was further improved by Lu, Wang and Yu
\cite{LWY20} to $n\geq 2kt$.
Keevash, Lifshitz, Long and Minzer \cite{KLLM1, KLLM2} independently verified this conjecture for $n > Ckt$ for some (large and unspecified) constant $C$.
Recently, Kupavskii \cite{K21} gave the concrete dependencies on the parameters by showing the conjecture holds for $n > 3ekt$ with $t > 10^7$.

For 3-graphs, Lu, Yu, and Yuan \cite{LYY2} proved
that, for sufficiently large $n$ with $n\equiv 0\pmod 3$,
if $\delta_1(F_i)>{n-1\choose 2}-{2n/3\choose 2}$ for $i\in [n/3]$
then ${\cal F}$ has a rainbow matching. This implies the result of  K\"uhn, Osthus, and Treglown  \cite{KOT13} and
Khan \cite{Kh16} on perfect matchings in 3-graphs.

In this paper, we prove the following result on rainbow matchings in
4-graphs,  which gives Khan's result \cite{Kh16} on perfect
matchings in 4-graphs as a special case.

\begin{theorem}\label{main}
Let $n$ be a sufficiently large integer with $n \equiv 0\pmod 4$. Let ${\cal F}=\{F_1,\ldots,F_{n/4}\}$ such that, for $i\in [n/4]$,
$F_i\subseteq {[n]\choose 4}$ and $\delta_1(F_i)>{n-1\choose
  3}-{3n/4\choose 3}$. Then ${\cal F}$ admits a rainbow matching.
\end{theorem}

The bound on $\delta_1(F_i)$ in Theorem \ref{main} is sharp.
To see this, let $k,m,n$ be positive integers, such that $k\ge
2$ and $2\leq m\leq n/k$. Let
%$H_k(n,m)\subseteq {[n]\choose k}$ su
\[
H_{k}(n,m)=\left\{e\in {[n]\choose k} : e\cap [m]\neq \emptyset
\mbox{ and } e\cap \left([n]\setminus [m]\right)\neq \emptyset
\right\}.
\]
and
\[
H_{k}^*(n,m)=\left\{e\in {[n]\choose k} : e\cap [m]\neq \emptyset
\right\}.
\]
Then $\delta_1(H_k(n,m))=\delta_1(H_k^*(n,m))={n-1\choose k-1}-{n-1-m\choose k-1}$ and
$H_k(n,m)$ has no matching of size $m+1$. It follows that when  $n\equiv
0 \pmod k$, we have $\delta_1(H_k(n,n/k-1))={n-1\choose k-1}-{n-n/k\choose k-1}$ and $\{H_k(n,n/k-1),\ldots,H_k(n,n/k-1)\}$ admits no rainbow matching.

  We prove Theorem~\ref{main} by working with a $5$-graph $H({\cal
    F})$  obtained
  from ${\cal F}=\{F_1,\ldots,F_{n/4}\}$: The vertex set of $H({\cal
    F})$ is $[n]\cup \{x_1, \ldots, x_{n/4}\}$ and the edge set of
  $H({\cal F})$ is  $\bigcup_{i=1}^{n/4}\{e\cup \{x_i\} :  e\in
  E(F_i)\}$. Clearly,  ${\cal F}$ admits a rainbow matching if and only if $H({\cal
    F})$ has a perfect matching.

For convenience,  we say that  a $(k+1)$-graph $H$ is \emph{$(1,k)$-partite} if there exists a
partition of $V(H)$ into sets $V_1, V_2$ (called {\it partition  classes}) such that for any $e\in E(H)$, $|e\cap V_1|=1$ and $|e\cap V_2|=k$.
A $(1,k)$-partite $(k+1)$-graph with partition classes $V_1,V_2$ is \emph{balanced}
if $k|V_1|=|V_2|$.  Thus, for instance,  $H({\cal F})$ above is a balanced
(1,4)-partite 5-graph with partition classes $X, [n]$.

More generally,
let $\mathcal{F}=\{F_1,\ldots,F_{m}\}$ be a family of $n$-vertex
$k$-graphs on a common
 vertex set $V$ and let $X=\{x_1, \ldots, x_m\}$ be a set disjoint
 from $V$. We use  $\mathcal{H}_{n,m}^k(\mathcal{F})$ to represent  the balanced
 $(1,k)$-partite $(k+1)$-graph with partition classes $X, V$
and edge set $\bigcup_{i=1}^m \{e\cup
\{x_i\} :  e\in E(F_i)\}$.
%% All the H_k(n,m) are on the same set of vertices.
If  $F_i=H_k(n,m)$ (or $H_k^*(n,m)$) for all $i\in [m]$,
then we write $\mathcal{H}_k(n,m)$ (or $\mathcal{H}_k^*(n,m)$) for
$\mathcal{H}_{n,m}^k(\mathcal{F})$ (or $\mathcal{H}_k^*(n,m)$). Now Theorem~\ref{main} is a direct consequence of the following result.

\begin{theorem}\label{general}
Let $n$ be an integer such that $n \equiv 0\pmod 4$ and $n$ is sufficiently large.
Let $H$ be a balanced
$(1,4)$-partite  5-graph with partition classes $X, [n]$ such that
 $\delta_1(N_H(x))>{n-1\choose
  3}-{3n/4\choose 3}$ for all $x \in X$.
Then $H$ admits a perfect matching.
\end{theorem}

Our proof of Theorem~\ref{general} is divided into two parts by
considering whether $H$ is close to $\mathcal{H}_4(n,n/4)$ or
not. For any real $\varepsilon>0$ and two $k$-graphs $H_1,H_2$ on the
same vertex set $V$, we say that $H_2$ is {\it $\varepsilon$-close} to $H_1$
if there exists an isomorphic copy $H_2'$ of $H_2$ with $V(H_2')=V$ such that $|E(H_1)\setminus E(H_2')|< \varepsilon |V(H_1)|^k$.

In Section 2, we prove Theorem~\ref{general} for the case
when $H$ is close to $\mathcal{H}_4(n,n/4)$. In fact,  we
 prove Theorem~\ref{general} for $(1,k)$-partite
$(k+1)$-graphs that are close to $\mathcal{H}_k(n,n/k)$, for all $k\ge 2$.
To prove Theorem~\ref{general} for the case
when $H$ is not close to $\mathcal{H}_4(n,n/4)$,
we will need to find a small
``absorbing'' matching in $H$, and this part is done in Section 3.
In Section 4, we show that if $H$ is not close to
$\mathcal{H}_4(n,n/4)$ then we can find a subgraph of $H$ that is
almost regular (in terms of vertex degree) and has maximum 2-degree
bounded above by $n^{0.1}$.  We make use of a recent stability
result of Gao, Lu, Ma, and Yu \cite{GLMY} for 3-graphs (see
Lemma~\ref{stable-lem}) and another result there on almost regular
spanning subgraphs.  We then complete the proof using a result of
Pippenger and Spencer.

\section{Hypergraphs close to $\mathcal{H}_k(n,n/k)$}

In this section, we prove Theorem~\ref{general}  for the case when $H_{n,n/k}^k(\mathcal{F})$
is $\varepsilon$-close to $\mathcal{H}_k(n,n/k)$ for some sufficiently small
$\varepsilon$.

We first prove Theorem~\ref{general} for those balanced $(1,k)$-partite
$(k+1)$-graphs $H$ in which, for each vertex $v\in V(H)$,  most edges of $H$ containing
$v$ also belong to $\mathcal{H}_k(n,n/k)$. More precisely, given
$\alpha>0$ and two $(k+1)$-graphs  $H_1,H_2$ on the same vertex set,  a
vertex $v\in V(H_1)$ is \emph{$\alpha$-bad} with respect to $H_2$
if $|N_{H_2}(v)\setminus N_{H_1}(v)|>\alpha |V(H_1)|^{k}$. (A vertex
$v\in V(H_1)$ is \emph{$\alpha$-good} with respect to $H_2$ if it is
not $\alpha$-bad with respect to $H_2$.) So if $v$ is $\alpha$-good
with respect to $H_2$ then all but at most
$\alpha |V(H_1)|^{k}$ of the edges containing $v$ in $H_2$ also lie
in $H_1$.

\begin{lemma}\label{good}
Let $k\ge 2$ be an integer,  $0<\alpha < (10^k k^k(k+1)!)^{-1}$, and let $n$ be an integer with $1/n\ll \alpha$ and $n\equiv 0\pmod k$. If
$H$ is a balanced $(1,k)$-partitie $(k+1)$-graph on the same vertex
set as $\mathcal{H}_{k}(n,n/k)$ and  every vertex of $H$ is
$\alpha$-good with respect to $\mathcal{H}_{k}(n,n/k)$, then $H$ has a perfect matching.
\end{lemma}

\pf Let $X, [n]$ denote the partition classes of $H$, and let
$W=[n/k]$ and $U=[n]\setminus W$.
Let $M$ denote a matching in $H$ such that $|e\cap X|=|e\cap W|=1$ for
each $e\in M$ and, subject to this, $|M|$ is maximum.
%consisting of edges of type $QWV^{(k-1)}$, where $V^{(k-1)}$ denotes
%$k-1$ $V$'s.
Let $U'=U\setminus V(M), W'=W\setminus V(M)$, and $X'=X\setminus
V(M)$. We may assume $|M|<n/k$; for otherwise, the assertion of the lemma is
true.

Note that $|M|\geq n/2k$.
For, suppose $|M|<n/2k$.
Then $|U'|/(k-1)=|W'|=|X'|\geq
n/2k$.
By the maximality of $|M|$, $H$ has no edge contained in $X' \cup W' \cup U'$.
Hence, for any $u\in U'$, we have
\begin{eqnarray*}
%\begin{split}
  & &|N_{\mathcal{H}_k(n,n/k)}(u)\setminus N_{H}(u)|\\
& \geq& |X'||W'|{|U'|\choose k-2} \\
& \ge & (n/2k)(n/2k)((k-1)n/2k-k+3)^{k-2}/(k-2)! \\
& \ge & \frac{n^k}{4k^25^{k-2}(k-2)!}  \quad (\mbox{since $n\ge 20k^2$
  and $k\ge 2$})\\
& >&\alpha \left(\frac{(k+1)n}{k}\right)^k=\alpha |V(H)|^k \quad
(\mbox{since  $\alpha < k^k/(4k^25^{k-2}(k-2)!(k+1)^k)$}).
%\end{split}
\end{eqnarray*}
Thus, $u$ is not $\alpha$-good with respect to $\mathcal{H}_k(n,n/k)$, a contradiction.

Fix $x\in X'$, $u_1,\ldots,u_{k-1}\in U'$, and $w\in W'$. Write $S=\{x,w,u_1,\ldots,u_{k-1}\}$.
If there exist distinct $e_1,\ldots, e_k \in M$ such that $H[S\cup (\cup_{i=1}^{k}e_i)]$
has a matching $M'$ of size $k+1$
such that for any $f \in M'$, $|f \cap X|=1=|f \cap W|$,
then $(M\setminus \{e_i: i\in
[k]\})\cup M'$ contradicts the choice of $M$.
So such $M'$ does not exist for any choice of distinct $e_1,\ldots, e_k \in M$. This implies that there exists a
$(k+1)$-subset $f$ of $V(H)$ such that $f\subseteq S\cup
(\bigcup_{i=1}^ke_i)$, $|f\cap
X'|=1=|f\cap W'|$, $|f\cap e_i|=1$ for $i\in [k]$, but $f\notin E(H)$.

Hence there exists $v\in S$ such that
\[
|N_{\mathcal{H}_k(n,n/k)}(v)\setminus
N_{H}(v)|
>\frac{1}{k+1}{n/2k\choose k}
>\frac{(n/2k-k+1)^k}{(k+1)!}
>\frac{(n/3k)^k}{(k+1)!}
>\alpha n^k,
%{n/2k\choose k}{n(1-1/k)/2\choose k(k-1)}>\alpha n^k,
\]
since $n > 6k(k-1)$ and $\alpha < (3^kk^k(k+1)!)^{-1}$. This is a contradiction.
\qed

\medskip

 To achieve the goal of this section, we need a result from \cite{LYY2}.

\begin{lemma}[Lu, Yu, and Yuan \cite{LYY2}]\label{LYY}
Let $n,t, k$ be positive integers such that $n> 2k^4 t$.
For $i\in [t]$, let $G_i\subseteq {[n]\choose k}$ such that $\delta_1(G_i)>{n-1\choose k-1}-{n-t \choose k-1}$.
Then $\{G_1, \ldots, G_t\}$ admits a rainbow matching.
\end{lemma}

%% We should be able to relax this result to $k \ge 2$.

\begin{lemma}\label{close}
Let $k\ge 3$ be an integer,  $0< \varepsilon <(10k)^{-6}$, and let $n$
be an integer with $n\equiv 0\pmod k$ and $n\ge 20k^2$. Let
$H$ be a balanced $(1,k)$-partitie $(k+1)$-graph with partition
classes $X, [n]$ and $V(H)=V(\mathcal{H}_{k}(n,n/k))$.  If $H$ is  $\varepsilon$-close to
$\mathcal{H}_k(n,n/k)$ and $\delta_1(N_{H}({x}))>{n-1\choose
  k-1 }-{n-n/k\choose k-1}$ for all  $x\in X$, then
$H$ has a perfect matching.

\end{lemma}

\pf
Let $W=[n/k]$ and $U=[n]\setminus [n/k]$ be the partition classes of
$H_k(n, [n/k])$ in $\mathcal{H}_{k}(n,n/k)$.
Let $B$ denote the set of $\sqrt{\varepsilon}$-bad vertices in $H$
with respect to $\mathcal{H}_{k}(n,n/k)$.
Since $H$ is $\varepsilon$-close to $\mathcal{H}_{k}(n,n/k)$, we have
$|B|\leq 2(k+1)\sqrt{\varepsilon}n$;
otherwise,
\begin{align*}
|E(\mathcal{H}_{k}(n,n/k))\setminus E(H)|&\geq \frac{1}{k+1}\sum_{v\in V(H)}|N_{\mathcal{H}_{k}(n,n/k)}(v)\setminus N_H(v)|\\
&> 2(k+1) \sqrt{\varepsilon}  n\cdot \frac {1}{k+1} \sqrt{\varepsilon}|V(H)|^{k}\geq \varepsilon |V(H)|^{k+1},
\end{align*}
a contradiction.

Let $U^{b}=U\cap B$, $X^{b}=X\cap B$, and $W^{b}=W\cap B$. Let
$W^{g}=W\setminus W^{b}$.  For convenience, write  $q=|X^{b}|$, $w^{b}=|W^{b}|$, and
$r=q+w^{b}$.
Moreover, let $x_1, \ldots, x_r$ be distinct such that $X^{b}=\{x_{1},\ldots,x_{q}\}$,
let $W' \subseteq W^g$ be a set of size $n/k-r$, and  let $G_i=N_{H}(x_i)-W^{'}$ for
$i\in [r]$. Then, for $i\in [r]$,
\begin{align*}
\delta_1(G_i)&\geq \delta_1(N_{H}(x_i))-\left({n-1\choose k-1}-{n-|W^{'}|-1\choose k-1}\right)\\
&>{n-|W^{'}|-1\choose k-1} -{n-n/k\choose k-1} \\
&={n-n/k+r-1\choose k-1} -{n-n/k\choose k-1}.
\end{align*}

Thus, by Lemma \ref{LYY} (with $n-n/k+r$ as $n$ and $r$ as $t$), $\{G_{1},\ldots,G_{r}\}$ admits a rainbow
matching, say $M_0=\{e_i\in E(G_i): i\in [r]\}$. Now $M_0'=\{e_i\cup \{x_i\}:  i\in [r]\}$ is a matching in $H$ covering $X^b$.
(Note this is the only place in this proof that requires the degree condition in the statement.)

Let $H_1=H-V(M_0')$.
Since $r\leq |B| \le 2(k+1)\sqrt{\varepsilon}n$ and $\varepsilon < (10k)^{-6}$, every vertex in $X\setminus V(M_0')$ is
$\varepsilon^{1/3}$-good with respect to $\mathcal{H}_{k}(n,n/k)-V(M_0')$. %Let $B_1=B-V(M_0')$.
Choose $\eta$ such that $0<\varepsilon\ll \eta\ll 1/k$, and let
$$B':=\{v\in B\setminus V(M_0'): |\{e\in E(H): v\in e \mbox{ and }
|e\cap W^g|=1\}|\ge \eta n^k\}.$$
Since $|B'|\leq |B|\leq 2(k+1)\sqrt{\varepsilon}n$ and $\varepsilon\ll
\eta$, we may greedily pick a matching $M_1$ in $H-V(M_0')$ such that
$B'\subseteq V(M_1)$ and every edge in $M_1$
contains at least one vertex from $B'$ and exactly one vertex from $W^{g}$.

Now consider $H_2=H_1-V(M_1)$. Note that since $n > 20k^2$,
\begin{equation*}
\begin{split}
\delta_1(H_2)
&\geq \delta_1(H)-(k+1)|M_0'\cup M_1|n^{k-1} \\
&> \frac{n}{k}\left({n-1\choose k-1}-{n-n/k\choose k-1}\right)-2(k+1)^2\sqrt{\varepsilon}n^k.
\end{split}
\end{equation*}
Thus, for each $v\in B\setminus V(M_0'\cup M_1)$,  the number of edges
of $H_2$ containing $v$ and no vertex of $W^{g}$ is at least
%(since $n > {k}/{(1-0.9^{1/(k-1)})}$),
\begin{eqnarray*}
%\begin{split}
& &\delta_1(H_2)-\eta n^k - \sum_{i=2}^{k-1}\frac{n}{k} {|W^g| \choose i} {n - |W^g|-1 -k|M_0'\cup M_1| \choose k-1-i}\\
&\geq & \frac{n}{k}\left({n-1\choose k-1}-{n-n/k\choose k-1}\right)-2(k+1)^2\sqrt{\varepsilon}n^k -\eta n^k - \frac{n}{k}\sum_{i=2}^{k-1}{n/k \choose i} {n - n/k-1 \choose k-1-i}\\
&\geq & \frac{n}{k}\left({n-1\choose k-1}-{n-n/k\choose k-1}\right)-2(k+1)^2\sqrt{\varepsilon}n^k -\eta n^k \\
& &- \frac{n}{k} \left( {n-1 \choose k-1} - {n-n/k-1 \choose k-1} - {n/k \choose 1} {n-n/k-1 \choose k-2} \right)\\
%&=\frac{n}{k}\left({n-n/k\choose k-1}+{n-n/k-1\choose k-1}-{n-n/k\choose k-1}\right)-3(k+1)^2\sqrt{\varepsilon}n^k -\eta n^k \\
&=& \frac{n}{k}\left(\frac{n}{k} - 1\right){n-n/k - 1\choose k-2}-2(k+1)^2\sqrt{\varepsilon}n^k -\eta n^k \\
&>& \eta n^k \quad \mbox{(since $\varepsilon\ll \eta\ll 1/k$)}.
%\end{split}
\end{eqnarray*}

Hence, we may greedily pick a matching $M_2$ in $H_2$ such that every edge in
$M_2$ contains at least one vertex from $B\setminus V(M_0'\cup M_1)$ and no vertex from $W^{g}$.

It is easy to see that $|M_0'\cup M_1\cup M_2|\leq 2(k+1)\sqrt{\varepsilon}n$.
Hence, every vertex of $H_2-V(M_2)$ is $\varepsilon^{1/4}$-good with
respect to $\mathcal{H}_{k}(n,n/k)-V(M_0'\cup M_1\cup M_2)$. Thus for
every vertex $u\in U\setminus V(M_0'\cup M_1\cup M_2)$, the number
of edges containing $u$ and exactly two vertices of $W\setminus
V(M_0'\cup M_1\cup M_2)$ as well as avoiding $V(M_0'\cup M_1\cup M_2)$  is at least
\[
\frac{n}{k}{n/k\choose 2}{n-n/k-1\choose k-3}-\varepsilon^{1/4}\left(n+ \frac{n}{k} \right)^k-(k+1)|M_0'\cup M_1\cup M_2|n^{k-1} >\eta n^k.
\]
Thus we may greedily pick a matching $M_2'$ such that $|M_2'|=|M_2|$ and every edge  of $M_2'$ contains exactly two vertices from $W^{g}$.

Put $M:=M_0'\cup M_1\cup M_2\cup M_2'$ and $m:=|M|$.  Let $H_3:=H-V(M)=H_2-V(M_2\cup M_2')$.
One can see that every vertex of $H_3$ is $\varepsilon^{1/5}$-good
with respect to $\mathcal{H}_k(n-km,n/k-m)=\mathcal{H}_k(n,m)-V(M)$. By Lemma \ref{good},
$H_3$ contains a perfect matching, say $M_3$.  Now $M_3\cup M$ is a perfect
matching in $H$.\qed

\section{Absorbing matching}

To deal with balanced $(1,k)$-partite $(k+1)$-graphs that are not
close to $\mathcal{H}_k(n, n/k)$, we need to find a small matching
that can ``absorb'' small sets of vertices. To find such a matching, we
need to use Chernoff inequality to bound deviations, see \cite{AS08}.

 \begin{lemma}[Chernoff inequality for small deviation]\label{chernoff1}
 Let $X=\sum_{i=1}^n X_i$, where each random variable $X_i$ has
 Bernoulli distribution with expectation $p_i$. For $\alpha \le 3/2$,
 $$ \mathbb{P}(|X-\mathbb{E}X| \ge \alpha \mathbb{E}X ) \le 2e^{-\frac{\alpha^2}{3}\mathbb{E}X}.$$
 In particular, when $X \sim Bi(n,p)$ and $\lambda < \frac{3}{2}np$, then
 $$ \mathbb{P}(|X-np| \ge \lambda ) \le e^{-\Omega(\lambda^2/(np))}.$$
 \end{lemma}

We now prove a $(1,k)$-partite version of the absorption lemma for
$(1,3)$-partite 4-graphs proved in \cite{LYY2}. Our proof follows along the same lines as in \cite{LYY2}.

\begin{lemma}\label{absorbing}
%$\zeta > 0$ be a real number an
	Let $k\geq 3$ be an integer and $0< b < 1/k$ be a
        constant. There exists an integers $n_1=n_1(k,b)$ such that
        the following holds for  any integer $n\ge n_1$:
Let $H$ be a $(1,k)$-partite $(k+1)$-graph with partition classes $X, [n]$ such that
$|X|=n/k$ and $\delta_1(N_H(x))>(1/2+b){n-1 \choose k-1 }$ for $x\in X$.
Then for any $c$ satisfying
$0<c<  \min \{( \frac{k^kc n}{6b^k (k!)^{k}})^{2},(2k^3(k+1)c^2)^{-10}\}$,
%$0<c<\zeta^{2k} (12k^22^k (k!)^k)^{-2}$,
there exists  a matching $M$ in $H$ %$\mathcal{H}$
	such that $|M|\le 2k c n$ and, for any balanced subset $S\subseteq V(H)$ %$S\subseteq V(\mathcal{H})$
       with $|S|\le (k+1)c^{1.5} n/2$,       $H[V(M)\cup S]$ %$\mathcal{H}[V(M)\cup S]$
         has a perfect matching.
\end{lemma}

\pf For balanced $R\in \binom{V(H)}{k+1}$ and balanced $Q\in \binom{V(H)}{k(k+1)}$,
	we say that $Q$ is \emph{$R$-absorbing}
	if both  $H[Q]$ and $H[Q\cup R]$ have perfect
        matchings. For each balanced $R\in \binom{V(H)}{k+1}$, let $\mathcal{L}(R)$ denote the collection of all $R$-absorbing sets in $H$.

\medskip
\textbf{Claim 1.~} For each balanced  $R\in \binom{V(H)}{k+1}$,
$|{\cal L}(R)|\ge b^k{n\choose k}^{k+1}/(2 (k^2)!)$.
%the number of  $R$-absorbing sets in $H$ is at least ${ b{n\choose k}^{k+1}}/{(2k 30^k (k^2)!)}$.

Let $R\in \binom{V(H)}{k+1}$ be a fixed balanced set, and let $R\cap
X=\{x\}$.  Note that the number of edges in $H$ containing
$x$ and intersecting $R\setminus \{x\}$ is at most $k{n\choose k-1}$.
Thus, since $\delta_1(N_H(x))>(1/2+b){n-1 \choose k-1 }$, the number of
edges  $e\in E(H)$ with $x\in e$ and $e\cap R=\{x\}$  is at least
\begin{align*}
\frac{n(1/2+b) {n-1 \choose k-1 }}{k} - k{n\choose k-1}\geq \frac{1}{2}{n\choose k}.
\end{align*}

%Let $R=\{x,u_1,\ldots,u_k\}$ be fixed with $x\in X$ and $u_i\in [n]$ for $i\in [k]$. Note that the number of edges in $H$ containing
%$x$ and intersecting $\{u_1, \ldots, u_k\}$ is at most $k{n\choose k-1}$,
%and $\delta_1(N_H(x_i))>(1/2+b){n-1 \choose k-1 }$.
%So the number of edges  $\{x,v_1,\ldots, v_k\}$ in $H$ such that
%$v_i\in [n]$ for $i\in [k]$ and $\{v_1, \ldots, v_k\}\cap \{u_1,\ldots, u_k\}=\emptyset$ is at least
%\begin{align*}
%\frac{n}{k}(\frac{1}{2}+b) {n-1 \choose k-1 } - k{n\choose k-1}\geq \frac{1}{2}{n\choose k}.
%\end{align*}
%Fix a choice of an edge $\{x,v_1,\ldots, v_k\}$ in $H$ such that
%$v_i\in [n]$ for $i\in [k]$ and $\{v_1, \ldots, v_k\}\cap \{u_1,\ldots, u_k\}=\emptyset$, and
%let $W_0 = \{v_1,\ldots, v_k\}$.

Fix a choice of $e\in E(H)$ with $x\in e$ and $e\cap R=\{x\}$, and write $R\setminus \{x\}=\{u_1, \ldots, u_k\}$
and $e\setminus \{x\}=\{v_1, \ldots, v_k\}$.   Let
$W_0=e\setminus \{x\}$. For each pair $\{u_j,v_j\}$, we choose a $k$-set $U_j$ disjoint from
$W_{j-1}\cup R$ such that both $U_j\cup \{u_j\}$ and $U_j\cup\{v_j\}$ are edges in $H$, and let $W_{j} :=U_{j}\cup W_{j-1}$.
If $W_k$ is defined then  $W_k$ gives an absorbing $k(k+1)$-set for $R$.

Note that for $j\in [k]$ there are $k+1 +jk$ vertices in $W_{j-1} \cup R$.  Thus, the number of edges in $H$ containing
 $u_j$ (respectively, $v_j$) and another vertex in $W_{j-1} \cup R$ is at most $(k+1 +jk){n\choose k-2} \frac{n}{k} <(k+1)n{n\choose k-2}$.
Since $\delta_1(N_H(x))>(1/2+b){n-1 \choose k-1 }$ for $x\in X$,
the number of sets $U_j$ for which $U_j \cap (W_{j-1} \cup R) = \emptyset$ and  both $U_j\cup \{u_j\}$ and $U_j\cup \{v_j\}$ are edges in $H$ is at least

\begin{eqnarray*}
%\begin{split}
& &\frac{n}{k} \left(2 (1/2+b){n-1 \choose k-1 }   - {n-1 \choose k-1 } \right) - 2(k+1)n{n \choose k-2} \\
& =&2b  {n \choose k }  - 2(k+1)n{n \choose k-2} \\
& >& b {n \choose k}
%\end{split}
\end{eqnarray*}
because $n$ is sufficiently large.
%$${n-1 \choose k-1 } / {n-t\choose k-1 } > (\frac{n-k+1}{n-bn/k})^{k-1} > (\frac{n-bn/k^3}{n-bn/k})^{k-1} = (1+\frac{1}{k} + \frac{1}{k^2})^{k-1} \ge 169/81.$$

To summarize, the number of $W_k$ defined above
from $e$ is at least $\left(b{n\choose k}\right)^k$. Hence,
there are at least  $\frac{1}{2}{n\choose k}(b{n\choose k})^k$
absorbing, ordered $k(k+1)$-sets for $R$,  with at most $(k^2)!$ of
them corresponding to a single $R$-absorbing set. Therefore,
\[
\mathcal{L}(R)\geq \frac{\frac{1}{2}{n\choose k}(b{n\choose k})^k}{(k^2)!}= \frac{b^k {n\choose k}^{k+1}}{2 (k^2)!}.
\]
This completes the proof of Claim 1.

\medskip

Now, let $c$ be a fixed constant with $0<c< \min \{( \frac{b^kk^k}{6(k!)^{k}})^{2},({2k^3(k+1)c^2})^{-10}\}$, and
choose a family ${\cal G}$ of balanced $k(k+1)$-sets of $V(H)$ by selecting each
of the ${ n/k \choose k}{n\choose k^2}$ balanced $k(k+1)$-subsets of
$V(H)$ independently with probability
\[
p := \frac{c n}{{n/k \choose k}{n\choose k^2}}.
\]
Then $\mathbb{E}(|{\cal G}|)=cn$ and $\mathbb{E}(|{\cal L}(R)\cap
{\cal G}|)=p|{\cal L}(R)|\ge b^k{n\choose k}^{k+1}cn/\left(2(k^2!){n/k\choose k}{n\choose
  k^2}\right)$.
It follows from Lemma~\ref{chernoff1} that, with probability $1-o(1)$,
\begin{align}\label{absor-eq1}
|{\cal G}|\leq 2c n
\end{align}
and, for all balanced $(k+1)$-sets $R$,
\begin{align}\label{absor-eq2}
|\mathcal{L}(R)\cap {\cal G}|\geq p |\mathcal{L}(R)|/2\geq \frac{b^kk^kc n}{6(k!)^{k}} \geq c^{1.5} n.
\end{align}
 Furthermore, the expected number of intersecting pairs of
 $k(k+1)$-sets in ${\cal G}$ is at most
\begin{align*}
{ n/k \choose k}{n\choose k^2}k(k+1)\left({ n/k -1\choose k-1}{n\choose k^2}+{n/k \choose k}{n-1\choose k^2-1}\right)p^2 \le 2k^3(k+1)c^2n\leq c^{1.9} n.
\end{align*}
Thus, using Markov's inequality, we derive that with probability at least $1/2$
\begin{align}\label{absor-eq3}
\mbox{${\cal G}$ contains at most $2c^{1.9} n$ intersecting pairs of $k(k+1)$-sets. }
\end{align}

Hence, there exists a family ${\cal G}$ satisfying (\ref{absor-eq1}),
(\ref{absor-eq2}), and (\ref{absor-eq3}). Delete one $k(k+1)$-set from
each intersecting pair in such a family ${\cal G}$, and remove all
non-absorbing $k(k+1)$-sets from ${\cal G}$. The resulting family,
call it ${\cal G}'$, consists of pairwise disjoint balanced, absorbing
$k(k+1)$-sets, and  satisfies
\[
|\mathcal{L}(R)\cap {\cal G}'|\geq c^{1.5} n/2,
\]
for all balanced $(k+1)$-sets $R$.

Since ${\cal G}'$ consists only of absorbing $k(k+1)$-sets, $H[V
({\cal G}')]$ has a perfect matching, say  $M$. By  (\ref{absor-eq1}),
$|M|\le 2kcn$.  For a balanced set $S\subseteq V(H)$ of size $|S|\leq
(k+1)c^{1.5} n/2$, we partition  $S$ into balanced $(k+1)$-sets $R_1,
\ldots, R_t$, where $t\le c^{1.5} n/2$. Since  $|\mathcal{L}(R_i)\cap {\cal G}'|\geq c^{1.5} n/2$,
there are distinct absorbing $k(k+1)$-set $Q_1, \ldots, Q_t$ in ${\cal
  G}'$ such that $Q_i$ is an $R_i$-absorbing set for $i\in [t]$. Now
$\mathcal{H}[V(M)\cup S]$ has a perfect matching which consists of a
perfect matching from each $H[Q_i\cup R_i]$.
\qed

\section{Fractional Perfect Matchings}
%\section{Hypergraphs not close to  $\mathcal{H}_4(n,n/4)$}

To deal with hypergraphs that are not close to
$\mathcal{H}_4(n,n/4)$, we need to control the independence number of
those hypergraphs. This is done in the same way as in \cite{LYY2} by
applying the hypergraph container result in  \cite{BMS15, ST15}.

First, we need the following lemma,  which is more general and
slightly stronger than Lemma 4.2 in \cite{LYY2} but with very similar
proof. Let $H$ be a hypergraph, $\lambda>0$ be a real number, and
${\cal A}$ be a family of subsets of $V(H)$. We say that $H$ is
\textit{$(\mathcal {A}, \lambda)$-dense} if $e(H[A])\ge \lambda e(H)$
for every $A \in \mathcal{A}$.

% We actually require the extremal graph to exclude the edges only in [t]??
\begin{lemma}\label{dense}
Let $\varepsilon$ be a constant such that
%$n \le 3k^2t$,
$0<\varepsilon \ll 1$, and let $n,k$ be integers such that $k \ge 3$ and $n\geq 40k^2/\varepsilon$. Let $a_1=\varepsilon/8k, a_2=\varepsilon/8k^3$,
and $a_3 < \varepsilon/(2^k \cdot k! \cdot 8k)$.
Let $H$ be a $(1,k)$-partite $(k+1)$-graph with vertex partition classes  $X,[n]$ with $|X|=n/k$.
Suppose $d_H(\{x,v\})\geq {n-1 \choose k-1}-{n-n/k \choose k-1}- a_3
n^{k-1}$ for any $x\in X$ and $v \in [n]$,
and $|E(\mathcal{H}_k(n,n/k))\setminus E(H_0)|\geq \varepsilon
e(\mathcal{H}_k(n,n/k))$ for any isomorphic copy $H_0$ of $H$ with
$V(H_0)=V(\mathcal{H}_k(n,n/k))$. Then $H$
%If $H$ is not $\varepsilon$-close to $\mathcal{H}_k(n,n/k)$, then $H$
is $(\mathcal{A}, a_1)$-dense, where $\mathcal{A}=\{A\subseteq V(H) :
|A\cap X|\ge (1/k-a_1) n \mbox{ and }\ |A\cap [n]|\ge (1-1/k-a_2) n\}$.
\end{lemma}

\pf We prove this by way of contradiction. Suppose that there exists $A\subseteq V(H)$ such that $|A\cap X|\ge (1/k-a_1) n$,
$|A\cap [n]|\ge  (1-1/k-a_2) n$, and $e(H[A])\le a_1 e(H)$. Without loss of generality, we may choose $A$ such that $|A\cap X|= (1/k-a_1) n$ and
$|A\cap[n]|= (1-1/k-a_2) n$. Let $U\subseteq [n]$ such that $A \cap
[n] \subseteq U$ and $|U|=n-n/k$. Let $A_1=A\cap X$, $A_2=X\setminus
A$, $B_1=A \cap [n]$, and $B_2=U\setminus A$. Relabel the vertices of
$H$ in $[n]$ if necessary, so that $U=[n]\setminus [n/k]$.

 Let $H_0$ denote the isomorphic copy of $H$ with the same partition classes
 $X, [n]$ as $\mathcal{H}_k(n,n/k)$. % and $U = [n] \setminus  [n/k]$.
 We derive a contradiction by showing
that $|E(\mathcal{H}_k(n,n/k))\setminus E(H_0)|< \varepsilon e(\mathcal{H}_k(n,n/k))$.
Note that %since $n\le 3k^2t$,
\begin{align*}
e(\mathcal{H}_k(n,n/k)) = \frac{n}{k} \left( {n \choose k} - {n-n/k \choose k}- {n/k\choose k}\right)
 \ge \frac{n}{k} {n \choose k}/k,
 \end{align*}
and, further,
\begin{align}\label{lem41-eq1}
e(\mathcal{H}_k(n,n/k)) \ge \frac{n}{k} {n \choose k}/k =
\frac{n^2}{k^3}{ n-1 \choose k-1} > \frac{n^2}{k^3}{ n-n/k \choose
  k-1}>\frac{n^3}{k^4}{ n-n/k \choose k-2}.
  \end{align}
Also, since $n > 2k$,
\begin{align}\label{lem41-eq2}
e(\mathcal{H}_k(n,n/k)) \ge \frac{n}{k} {n \choose k}/k > \frac{n^{k+1} }{2^k \cdot k! \cdot k^2 }.
 \end{align}

Consider $x \in A_1$ and $v \in [n] \setminus [n/k]$.
Let $E_{H_0}(B_1,\{x,v\}) =\{e\in E(H_0): \{x,v\}\subseteq e\subseteq
B_1\cup \{x,v\}\}$. %denote the set of edges contained entirely in $B_1 \cup \{x, v\}$ in $H_0$ that must contain $x$ and $v$.
%|\{e | x \in e, e \in E(H_0), V(e) \subseteq B_1 \cup \{x\} \}|
%The number of edges in $H_0$ containing $x,v$ that also exist in $\mathcal{H}_k(n,n/k)$ is the number of edges in $H_0$ containing $x,v$ and
 %intersecting $[n/k]$.
Note that %$$\{e\in E(H_0): \{x,v\}\subseteq e\}\cap E({\cal H}_k(n,
%n/k)) \subseteq \{e\in E(H_0): \{x,v\}\subseteq e \mbox{ and } e\cap [n/k]\ne \emptyset\}.$$
% Hence,
 for $v \in B_1$, we have
\begin{eqnarray*}
 & &  |\{e\in E(H_0) : \{x,v\} \subseteq e \mbox{ and } e \cap [n/k] \neq \emptyset\}| \\
&\ge &  d_{H_0}(\{x,v\}) - |\{e \in E(H_0-[n/k]) : \{x,v\} \subseteq e
       \mbox{ and } e \cap B_2 \neq \emptyset\}| - |E_{H_0}(B_1,\{x,v\})| \\
&\ge &  \left( {n-1 \choose k-1}-{n-n/k \choose k-1}- a_3 n^{k-1} \right) - a_2 n {n-n/k \choose k-2} - |E_{H_0}(B_1,\{x,v\})|\label{lem41-eq3}.
\end{eqnarray*}
For $v \in B_2$, we have
\begin{eqnarray*}
 & &  |\{e\in E(H_0) : \{x,v\} \subseteq e \mbox{ and } e \cap [n/k] \neq \emptyset\}| \\
&\ge &  d_{H_0}(\{x,v\}) - |\{e \in E(H_0-[n/k]) : \{x,v\} \subseteq e \}|  \\
&\ge &  \left( {n-1 \choose k-1}-{n-n/k \choose k-1}- a_3 n^{k-1} \right) - {n-n/k \choose k-1} \label{lem41-eq4}.
\end{eqnarray*}

%By assumption, $d_{H_0}(x)\geq {n\choose k}-{n-t+1 \choose k}-c^{1/4} n^k$.
So we have
\begin{eqnarray}
 & &   \sum_{x \in A_1} \sum_{v \in [n] \setminus [n/k]} |\{e \in E(H_0),:
     \{x,v\} \subseteq e \mbox{ and } e \cap [n/k] \neq \emptyset\}| \nonumber\\
&\ge &  \sum_{x \in A_1} \sum_{v \in [n] \setminus [n/k]}\left( {n-1 \choose k-1}-{n-n/k \choose k-1}- a_3 n^{k-1}  - a_2 n {n-n/k \choose k-2} \right)\nonumber\\
 & & -\sum_{x \in A_1} \sum_{v \in B_2} {n-n/k \choose k-1} - \sum_{x \in A_1} \sum_{v \in B_1}|E_{H_0}(B_1,\{x,v\})| \nonumber\\
&=&  \sum_{x \in A_1} \sum_{v \in [n] \setminus [n/k]}\left( {n-1 \choose k-1}-{n-n/k \choose k-1}- a_3 n^{k-1}  - a_2 n {n-n/k \choose k-2} \right) \nonumber\\
 & & -|A_1||B_2|{n-n/k \choose k-1} - |E(H_0[A])|.\quad
\end{eqnarray}
Note that
\begin{align*}
a_2\frac{n^3}{k^2}  {n-n/k \choose k-2} -  |A_1||B_2|{n-n/k \choose k-1}&>a_2\frac{n^2}{k}\left(\frac{n}{k}  {n-n/k \choose k-2}-{n-n/k \choose k-1}\right)\\
&>a_2\frac{n^2}{k}\left(\frac{n}{k}  {n-n/k-1 \choose k-2}-{n-n/k \choose k-1}\right)\\
&=0,
\end{align*}
That is,
\begin{align}\label{eq-A1B2}
a_2\frac{n^3}{k^2}  {n-n/k \choose k-2} -  |A_1||B_2|{n-n/k \choose k-1}>0.
\end{align}
Therefore, we have

\begin{eqnarray*}
& & |E(\mathcal{H}_k(n,n/k))\setminus E(H_0)| \\
&= & \sum_{x \in A_1} |\{ e \in E(\mathcal{H}_k(n,n/k))\setminus E(H_0): x \in e\}| + \sum_{x \in A_2} |\{ e \in E(\mathcal{H}_k(n,n/k))\setminus E(H_0) : x \in e\}| \\
&\leq & \sum_{x \in A_1} \sum_{v \in [n] \setminus [n/k]} |\{e \in
        E(\mathcal{H}_k(n,n/k))\setminus E(H_0): \{x,v\} \subseteq
        e\}| + |A_2|\cdot e({ H}_k(n,n/k)) \\
&\le & \sum_{x \in A_1} \sum_{v \in [n] \setminus [n/k]} \left( {n-1
       \choose k-1} - {n-n/k \choose k-1} - |\{ e \in E(H_0) : \{x,v\} \subseteq e
       \mbox{ and }e \cap [n/k] \neq \emptyset\}|\right) \\
& &+ |A_2| \cdot e({ H}_k(n,n/k))  \\
&\le & \sum_{x \in A_1} \sum_{v \in [n] \setminus [n/k]}\left(a_3
       n^{k-1}  + a_2 n {n-n/k \choose k-2}\right)+|A_1||B_2|{n-n/k \choose k-1}+ |E(H_0[A])| \\
& &+k a_1 \cdot e({\cal H}_k(n,n/k)) \quad \mbox{(by (\ref{lem41-eq3}))}\\
%&\le & \sum_{x \in A_1} \sum_{v \in [n] \setminus [n/k]} \Big[ \left({n-1 \choose k-1} - {n-n/k \choose k-1}\right) - \Big( \Big( {n-1 \choose k-1}-{n-n/k \choose k-1}- a_3 n^{k-1} \Big) \\
%& &  - |B_2| {n-n/k \choose k-2} - |E_{H_0}(B_1,\{x,v\})| \Big)\Big] + a_1 n \cdot e(\mathcal{H}_k(n,n/k))/(n/k) \\
%&\le &\sum_{x \in A_1} \sum_{v \in [n] \setminus [n/k]} \left[  a_3 n^{k-1} + a_2 n {n-n/k \choose k-2} +  |E_{H_0}(B_1,\{x,v\})| \right] + (k a_1) \cdot e(\mathcal{H}_k(n,n/k)) \\
&\le & \frac{n}{k}\left(n-\frac{n}{k}\right)\left(a_3 n^{k-1}  + a_2 n
       {n-n/k \choose k-2}\right)+|A_1||B_2|{n-n/k \choose k-1}+a_1 e(H_0)+k a_1 \cdot e({\cal H}_k(n,n/k))\\
&< & \frac{n^2}{k}\left(a_3 n^{k-1}  + a_2 n
       {n-n/k \choose k-2}\right)+a_1 e(H_0)+ ka_1 \cdot e({\cal H}_k(n,n/k))\quad \mbox{(by (\ref{eq-A1B2}))}\\
       %+\frac{n}{k}a_2n{n-n/k \choose k-1}
%&= & |A_1| a_3  n^k  + |A_1| a_2 n^2 {n-n/k \choose k-2} + \sum_{x \in A_1} \sum_{v \in [n] \setminus [n/k]}  |E_{H_0}(B_1,\{x,v\})| + (k a_1) \cdot e(\mathcal{H}_k(n,n/k)) \\
%&\le & (2^k \cdot k! \cdot k a_3) \cdot e(\mathcal{H}_k(n,n/k)) + (k^3 a_2) \cdot e(\mathcal{H}_k(n,n/k)) \\
% & & + e(H_0[A]) + (k a_1) \cdot e(\mathcal{H}_k(n,n/k)) \\
&< & a_1 e(H_0) + \left( 2^k \cdot k! \cdot k a_3 +  k^3 a_2 +  ka_1\right) \cdot e(\mathcal{H}_k(n,n/k)) \quad \mbox{(by (\ref{lem41-eq1}) and (\ref{lem41-eq2}))} \\
&\le & a_1 \frac{n}{k} {n \choose k} +\left(2^k \cdot k! \cdot k a_3  + k^3 a_2 + ka_1 \right) \cdot e(\mathcal{H}_k(n,n/k)) \\
&\le & \left(k a_1+ 2^k \cdot k! \cdot k a_3 + k^3 a_2 + k a_1 \right) \cdot e(\mathcal{H}_k(n,n/k)) \quad \mbox{(by (\ref{lem41-eq2}))}\\
&\leq & \varepsilon \cdot e(\mathcal{H}_k(n,n/k)),
\end{eqnarray*}
a contradiction. \qed  % since $H$ is not $\varepsilon$-close to $\mathcal{H}_k(n,n/k)$.  \qed

\medskip

To prove a fractional matching lemma, we need a recent result of
Gao, Lu, Ma, and Yu \cite{GLMY} on Erd\H{o}s' matching conjecture for stable
graphs.
For a hypergraph $H$, let $\nu(H)$ denote the size of maximum matching in $H$.
{\L}uczak and
Mieczkowska \cite{LM} proved that there exists a positive integer $n_1$ such that
for integers $m,n$ with $n\ge  n_1$ and $1 \le m \le  n/3$,
if $H$ is an $n$-vertex 3-graph with $e(H) > \max\{{n\choose
  3}-{n-m+1\choose 3}, {3m-1\choose 3}\}$ then $\nu(H) \ge m$. The
result of Gao, Lu, Ma, and Yu \cite{GLMY} may be viewed as a stability version of
this {\L}uczak-Mieczkowska result.

For sets $e=\{u_1, ...,u_k\} \subseteq [n]$ and $f=\{v_1,...,v_k\} \subseteq [n]$ with
$u_i < u_{i+1}$ and $v_i < v_{i+1}$ for $i \in [k-1]$, we write $e \le
f$ if $u_i\leq v_i$  for all $i \in [k]$.
A hypergraph $H$ with $V(H) = [n]$ and $E(H) \subseteq {[n] \choose
  k}$ is said to be \textit{stable} if, for any $e, f \in {[n] \choose k}$
with $e \le f$, $f \in E(H)$ implies $e \in E(H)$.
The following is a special case (when $m=3n/4$) of Lemma 4.2 in \cite{GLMY}.
Note that one of the extremal configurations of Lemma 4.2 in \cite{GLMY}, namely $\mathcal{D}(n,m,3)$ (the $3$-graph with vertex set $[n]$ and edge set ${\lfloor 3n/4 - 1 \rfloor \choose 3}$) does not occur here.

\begin{lemma}[Gao, Lu, Ma, and Yu]\label{stable-lem} For
  any $\eta>0$ there exists $n_0>0$ with the
  following properties: %Let $n$ be a sufficiently large and let
                        %$0<\eta\ll\varepsilon\ll 1$.
Let $n$ be an integer with $n\ge n_0$, and let
$H$ be a stable 3-graph on the vertex set $[n]$.
If $e(H)>{n\choose
  3}-{3n/4\choose 3}-\eta^4n^3$ and $\nu(H) < n/4$, then
$H$ is $\eta$-close to $H_3^*(n,n/4-1)$.
%$S(n,m,3) or D(n,m,3).
\end{lemma}

We will use perfect fractional matchings in a hypergraph $H$ not close to
${\cal H}_4(n,n/4-1)$ to obtain an almost regular subgraph of $H$. A
{\it fractional matching} in $H$ is a function $f: E(H)\to
\mathbb{R}^+$, where $\mathbb{R}^+$ is the set of non-negative reals, such that $\sum_{v\in e}f(e)\le 1$ for all $v\in
V(H)$, and it is {\it perfect} if $\sum_{v\in e}f(e)=1$ for all $v\in
V(H)$. We write $\nu_f(H) =\max\{\sum_{e\in E(H)}f(e): f \mbox{ is a
  fractional matching in } H\}$.

\begin{lemma}\label{frac-PM}
For any $\varepsilon >0$ there exists
$0<\rho\ll \varepsilon$ and $n_0=n_0(\varepsilon)$ such that, for
any integer $n$ with $n\ge n_0$ and $n\equiv 0\pmod 4$, the following
holds:
Let $H$ be a balanced $(1,4)$-partite $5$-graph  with partition
classes $X,[n]$, such that
%and on the same vertex set as $\mathcal{H}_{4}(n,n/4)$.
$d_H(\{x,v\})\geq {n-1\choose 3}-{3n/4\choose 3}-3\rho n^3$ for any $x
\in X$ and $v\in [n]$ and $H$ contains no independent set $S$  with
$|S\cap X|\geq n/4-\varepsilon n$ and $|S\cap [n]|\geq
3n/4-\varepsilon n$, then $H$ contains a fractional perfect matching.

\end{lemma}

\pf  Let $X=\{x_1, \ldots, x_{n/4}\}$ and let $\omega:V(H)\rightarrow
R^+$ be a minimum fractional vertex cover of $H$, that is,
$\sum_{v\in e}\omega(v)\ge 1$ for all $e\in E(H)$ and, subject to
this, $\omega(H):=\sum_{v\in V(H)}\omega(v)$ is minimum. We may assume that the vertices in
$X$ and $[n]$ are labeled such that
$\omega(x_1)\geq \cdots \geq \omega(x_{n/4})$ and $\omega(1)\geq
\cdots \geq \omega(n)$.

Let $H'$ be the $(1,4)$-partite 5-graph with vertex set $V(H)$
and edge set
\[
E(H)\cup \left\{e\in {V(H)\choose 5} : |e\cap X|=1 \ \mbox{and }\sum_{v\in e}w(v)\geq 1\right\}.
\]
Thus, by definition, $\omega$ is also a vertex cover of
$H'$. Moreover, since $E(H)\subseteq E(H')$, $\omega$ is also minimum
fractional vertex cover of $H'$. Hence, by linear programming duality,
we have
  $\nu_f(H)=\omega(H)=\omega(H')=\nu_f(H')$. Therefore, it suffices to
  show that $H'$ has a perfect matching. First, we show that $H'$
  is stable.

\begin{itemize}
  \item [(1)] Let $T_1=\{x_{i_1},i_2,i_3,i_4,i_5\}$ and
    $T_2=\{x_{j_1},j_2,j_3,j_4,j_5\}$ be balanced 5-element subsets of
    $V(H)$, with $x_{i_1},x_{j_1}\in X$ and $i_l\geq j_l$ for $l\in
    [5]$.  Then $T_2\in E(H')$ implies that $T_1\in E(H')$.
\end{itemize}
Since $i_l\geq j_l$ for $l \in [5]$, we have $\omega(x_{i_1})\geq \omega(x_{j_1})$
and $\omega(i_l)\geq \omega(j_l)$ for $2\leq l\leq 5$.  Note that $ \sum_{v\in
  T_2}\omega(v)\geq 1$ as $T_2\in E(H')$. Thus $\sum_{v\in T_1}\omega(v)\geq \sum_{v\in T_2}\omega(v)\geq
1$. Hence,  $T_1\in E(H')$ by the definition of $H'$, which  completes
the proof of (1).

\medskip

Let $G$ be the $3$-graph with vertex set $[n-1] $ and edge set
$N_{H}(\{x_{n/4}, [n]\})$.
%We choose $\eta$ such that $\rho\ll \eta\ll \varepsilon$.
We may assume that

\begin{itemize}
\item [(2)] $G$ has no matching of size $n/4$.
\end{itemize}
%\textbf{Case 1.} $G$ has a matching $M$ of size $n/4$.
For, suppose $G$ has a matching of size $n/4$, say
$M=\{e_1,\ldots,e_{n/4}\}$. Partition $V(H)\setminus V(M)$ into
2-sets $f_1,\ldots,f_{n/4}$, such that $|f_i\cap X|=1$ for all $i\in
[n/4]$. By (1), $M\subseteq N_{H}(f_i)$. Thus $M'=\{e_i\cup f_i : i\in
[n/4]\}$ is a perfect matching in $H$, completing the proof of (2).

\medskip

%\textbf{Case 2.} $\nu(G)<n/4$.

%We choose  $\eta$ such that $\eta\ll \varepsilon$.
%We further choose $\rho>0$ with  $\rho\ll \eta$ such that
Since $d_H(\{x,v\})\geq {n-1\choose 3}-{3n/4\choose 3}-3\rho n^3$ for any $x
\in X$ and $v\in [n]$, $e(G)>{n-1\choose 3}-{3n/4\choose 3}-3\rho
n^3$. Hence, by (2) and by Lemma
\ref{stable-lem},
 $G$ is $\eta$-close to
$H_3^*(n-1,n/4-1)$ with respect to the partition $[n-1]\setminus [n/4-1],
[n/4-1]$, where $\eta = (3\rho)^{1/4}$.  Let $Y=[n/4-\lceil \eta n\rceil]$. We claim that

\begin{itemize}
\item [(3)] for every $y\in Y$, $d_G(y)\geq {n-1\choose 2}-4\sqrt{\eta}
  n^2$.
\end{itemize}
For, otherwise,
since $G$ is stable, $d_G(z) < {n-1 \choose 2} - 4\sqrt{\eta}n^2$ for $z \in Z := [n/4 - 1] \backslash [n/4 - \lceil \eta n \rceil]$.
Hence,
%Then $G$ contains at least $n/4-4\sqrt{\eta} n$ vertices $x$ such that $d_G(x)\geq {n-1\choose 3}-\sqrt{\eta} n^3$, otherwise, we have
\[
|E(H_3^*(n-1,n/4-1))\setminus E(G)|\geq  \frac{1}{3} \sum_{z \in Z} \left( {n-1 \choose 2} - d_G(z) \right) \ge \frac{1}{3} (\sqrt{\eta}n - 1)4\sqrt{\eta}n^2  >\eta (n-1)^3,
\]
a contradiction which completes the proof of (3).
%Let $X$ denote the set of vertices $x$ with $d_G(x)\geq {n-1\choose 3}-\sqrt{\eta} n^3$. Without loss generality, we may assume that $|X|=n/4-\lceil\sqrt{\eta}n\rceil$.

\medskip

Since $H$ contains no independent set $S$  such that
$|S\cap X|\geq n/4-\varepsilon n$ and $|S\cap [n]|\geq
3n/4-\varepsilon n$,
we may greedily  find a matching $M_1$ of size $\lceil\sqrt{\eta}n\rceil$ in $H-Y$.

Next we  greedily construct a matching $M_2$ of size $|Y|$ in
$G-V(M_1)$ such that $|e\cap Y|=1$ for all $e\in M_2$.
 For $y\in Y$, note that
$$|\{e\in E(G): y\in e \mbox{ and } e\cap V(M_1)\ne \emptyset\}|\le
3|M_1| n<3\lceil\sqrt{\eta } n \rceil n \le 4 \lceil\sqrt{\eta}\rceil n^2.$$
By (3),
$d_G(1)-4\sqrt{\eta} n^2\geq {n-1\choose 2}-8\sqrt{\eta}
n^2>0;$ so there exists an edge $e_1\in E(G)\setminus V(M_1)$ such that
$|e_1\cap Y|=1$. Now suppose we have found a matching
$\{e_1,e_2,...,e_r\}$ in $G-V(M_1)$  such that
$|e_i\cap Y|=1$ for all $i\in [r]$. If $r=n/4-\lceil\sqrt{\eta}n\rceil$,
then $\{e_1,\ldots,e_r\}$ gives the desired matching $M_2$. So assume $r<|Y|$.
Write $G_r:=(G-V(M_1))-(\cup_{i=1}^r e_i)$. Let $v\in Y \backslash (V(M_1)\cup (\cup_{i=1}^r e_i))$. Note that $|[n] \backslash (Y\cup V(M_1)\cup (\cup_{i=1}^r e_i))|> n/4$. Since $d_G(v)>{n-1\choose 2}-4\sqrt{\eta}n^2$, the number of edges $e$ in $G$ with $v\in e$ and $e\backslash \{v\} \subseteq [n] \backslash (Y \cup V(M_1) \cup (\cup_{i=1}^r e_i))$ is at least
\[
{|[n] \backslash (Y \cup V(M_1) \cup (\cup_{i=1}^r e_i))|\choose 2}-4\sqrt{\eta}n^2>{n/4\choose 2}-4\sqrt{\eta}n^2>0,
\]
% Since $G$ is $\eta$-close to
%$H_3(n-1,n/4-1)$ with respect to partition $[n-1]\setminus [n/4-1],
%[n/4-1]$,
%\[
%|\{e\in E(G): e\cap Y=\{r+1\}\}|\ge {n-|Y|\choose 3}-\eta n^3,
%\]
%%then the number of edges $e$ with $e\cap Y=\{r+1\}$ in $G_r$ is at
%%least
%we have
%\[
%|\{e\in E(G_r): e\cap Y=\{r+1\}\}|\ge {n-1\choose
%  3}-\sqrt{\eta}n^3-|Y|n^2\ge {n-|Y|-2r\choose 3}-\eta n^3-3|M_1|n^2>0.
%\]
%%\[
%d_{G_r}(r+1)\geq {n-1\choose 3}-\sqrt{\eta} n^3-3|M_1|n^2-3(n/4-\lceil\sqrt{\eta}n\rceil){n\choose 2}>0,
%\]
So there exists an edge $e_{r+1}$ in $G_r$ such that $|e_{r+1}\cap
Y|=1$, contradicting the maximality of $r$.
%Continuing this process for $n/4-\lceil\sqrt{\eta}n\rceil$ steps, we
%obtain the desired matching $M_2$.

 Let $M_2=\{e_1,\ldots,e_{n/4-\lceil\sqrt{\eta}n\rceil}\}$.
%be a matching of size $n/4-\lceil\sqrt{\eta}n\rceil$ of graph
%$G-V(M_1)$.
Note $V(H)\setminus (V(M_1)\cup V(M_2))$ contains
$n/4-\lceil\sqrt{\eta}n\rceil$ vertex-disjoint 2-set, say
$f_1,\ldots,f_{n/4-\lceil\sqrt{\eta}n\rceil}$, such that $|f_i\cap
X|=1$ for $i\in [n/4-\lceil\sqrt{\eta}n\rceil]$. By (1),  $M_2\subseteq
N_{H'}(f_i)$ for $i\in [n/4-\lceil\sqrt{\eta}n\rceil]$. Write
$M_2'=\{f_i\cup e_i : i\in [n/4-\lceil\sqrt{\eta}n\rceil]\}$. Then $M_1\cup M_2'$ is a perfect matching in $H$. \qed

\medskip

\section{Random Rounding}

We need a result of
Lu, Yu, and Yuan \cite{LYY} on  the independence number of
a subgraph of a balanced  $(1,k)$-partite $(k+1)$-graph induced by a
random subset of vertices. It is stated for in \cite{LYY} for
$(1,3)$-partite 4-graphs, but the same proof (which uses the hypergraph
container result) also works for
$(1,k)$-partite $(k+1)$-graphs.

\begin{lemma}[Lu, Yu, and Yuan]\label{indep}
        Let $l, \varepsilon', \alpha_1,\alpha_2$ be positive reals, let $\alpha>0$ with  $\alpha \ll \min\{\alpha_1,\alpha_2\}$,
       let $k,n$ be positive integers, and let
       $H$ be a $(1,k)$-partite $(k+1)$-graph  with partition classes $Q,P$ such that $k|Q|=|P|=n$, $e(H)\ge ln^{k+1}$,  and $e(H[F])\ge
        \varepsilon' e(H)$ for all $F\subseteq
        V(H)$ with $|F\cap P|\ge \alpha_1 n$ and $|F\cap Q|\ge \alpha_2 n$.
         Let $R\subseteq V(H)$ be obtained  by taking each vertex of
           $H$ uniformly at random with probability $n^{-0.9}$.
        Then, with probability at least $1-n^{O(1)}e^{-\Omega (n^{0.1})}$, every independent set $J$ in $H[R]$
             satisfies $|J\cap P|\le (\alpha_1 +\alpha+o(1))n^{0.1}$ or $|J\cap Q|\le (\alpha_2 +\alpha+o(1))n^{0.1}$.
\end{lemma}

We also need Janson's inequality to provide an exponential upper bound for the lower tail of a sum of dependent zero-one random variables. See Theorem 8.7.2 in \cite{AS08}.

\begin{lemma}[Janson \cite{AS08}] \label{janson}
Let $\Gamma$ be a finite set and $p_i \in [0,1]$ be a real for $i \in \Gamma$.
Let $\Gamma_{p}$ be a random subset of $\Gamma$ such that the elements are chosen independently with $\mathbb{P}[i \in \Gamma_p] = p_i$ for $i \in \Gamma$.
Let $S$ be a family of subsets of $\Gamma$.
For every $A \in S$, let $I_A = 1$ if $A \subseteq \Gamma_p$ and $0$ otherwise.
Define $X = \sum_{A \in S} I_A$,
$\lambda = \mathbb{E}[X]$,
and
$\Delta = \frac{1}{2}\sum_{A \neq B} \sum_{A \cap B \neq \emptyset} \mathbb{E}[I_A I_B]$.
%and $\bar{\Delta} = \lambda + 2\Delta$.
Then, for $0 \le t \le \lambda$, we have
$$\mathbb{P}[X \le \lambda - t] \le \exp(-\frac{t^2}{2\lambda +4 {\Delta}}). $$
\end{lemma}

Now,  we use Chernoff bound and Janson's inequality to prove a result on several properties of  certain random   subgraphs.

 \begin{lemma}\label{lem1-5}
        Let $n, k$ be integers such that $n\ge k\geq 3$,
       let $H$ be a $(1,k)$-partite $(k+1)$-graph with partition classes $A,B$
       and  $k|A| = |B| = n$,
       and let $A_3\subseteq A$ and $A_4\subseteq B$ with $|A_i|=n^{0.99}$ for $i=3,4$.
	Take $n^{1.1}$ independent copies of $R$ and denote them by $R^i$, $1\le i\le n^{1.1}$, where $R$ is chosen from $V(H)$ by taking each vertex uniformly at random with probability $n^{-0.9}$ and then deleting $O(n^{0.06})$ vertices so that $|R|\in (k+1) \mathbb{Z}$ and $k|R\cap A|=|R\cap B| $.
        For each $S\subseteq V(H)$, let $Y_S:=|\{i: \ S\subseteq R^i\}|$.
            Then, with probability at least $1-o(1)$,  all of the following statements hold:
        \begin{itemize}
            \item [$(i)$] $Y_{\{v\}}=(1 \pm n^{-0.01}) n^{0.2}$  for all  $v\in V(H)$.  %\sim n^{0.2},
            \item [$(ii)$] $Y_{\{u,v\}}\le 2$ for all $\{u, v\} \subseteq V(H)$.
            \item [$(iii)$] $Y_e\le 1$ for all  $e \in E(H)$.
            \item [$(iv)$] For all $i= 1, \dots ,n ^{1.1}$, we have
              $|R_i\cap A| =(1/k\pm o(n^{-0.04}))n^{0.1}$ and  $|R_i\cap B| =(1\pm o(n^{-0.04}))n^{0.1}$,
              % $m\le n/k$ and
            \item [$(v)$] Suppose $\rho$ is a constant with $0<\rho <1$ such that
             $d_{H}(\{x,v\})\ge {n-1\choose
                   k-1}-{n-n/k\choose k-1}-\rho n^{k-1}$ for all $x\in A$ and $v \in B$.  Then
                   %there exists $\rho' > 0$ such that
                for $1 \le i \le n^{1.1}$, and for $x\in R_i \cap A$ and $v\in R_i \cap B$, we have
               $$d_{R_i}(\{x,v\})>  {|R_i \cap B| - 1\choose k-1}-{|R_i \cap B|-|R_i \cap B|/k \choose k-1}-3\rho |R_i \cap B|^{k-1},$$
           \item[$(vi)$] $|R_i\cap A_j|= |A_j|n^{-0.9}\pm n^{0.06}$ for $1 \le i \le n^{1.1}$ and $j\in \{3,4\}$.
            %\item[$(vii)$] for any idenpendent set $S$ of $R_i$, $|S\cap R_i\cap Q|< t(1-\varepsilon/2)$ and $|S\cap [n]\cap R_i|< (n-t)(1-\varepsilon)$
%with probability at most $n^{O(1)}e^{-\Omega (n^{0.1})}$.
        \end{itemize}
    \end{lemma}

\pf
For $1 \le i \le n^{1.1}$ and $j\in \{3,4\}$, $\mathbb{E}[|R_i \cap A|] = n^{0.1}/k$, $\mathbb{E}[|R_i \cap B|] = n^{0.1}$ and $\mathbb{E}[|R_i \cap A_j|] = n^{-0.9}|A_j|$. Recall the assumptions $|A_3|=|A_4|=n^{0.99}$.
By Lemma~\ref{chernoff1}, we have
 \begin{itemize}
 \item []  $\mathbb{P}\left(\left||R_i \cap A| - n^{0.1}/k\right| \ge n^{0.06} \right) \le e^{-\Omega(n^{0.02})}$,
 \item [] $ \mathbb{P}\left(\left||R_i \cap B| - n^{0.1}\right| \ge n^{0.06} \right) \le e^{-\Omega(n^{0.02})}$, and
  \item [] $\mathbb{P}\left(\left||R_i \cap A_j| - |A_j|n^{-0.9}\right| \ge n^{0.06} \right) \le e^{-\Omega(n^{0.03})}$.
\end{itemize}
 Hence, with probability at least $1-O(n^{1.1})e^{-\Omega(n^{0.02})}$, $(iv)$ and $(vi)$ hold.

 For every $v\in V(H)$, $\mathbb{E}[Y_{\{v\}}]=n^{1.1} \cdot n^{-0.9}= n^{0.2}$. By Lemma~\ref{chernoff1} again,
 $$ \mathbb{P}\left(\left||Y_{\{v\}}| - n^{0.2} \right| \ge n^{0.19} \right) \le e^{-\Omega(n^{0.18})}.$$
 Hence, with probability at least $1-O(n)e^{-\Omega(n^{0.18})}$, $(i)$ holds.

For positive integers $p, q$, let $Z_{p,q} = \left|S \in  {V(H) \choose p} : Y_S \ge q \right|$. Then
 $$\mathbb{E}\left[Z_{p,q}\right] \le  {n \choose p} {n^{1.1} \choose q} (n^{-0.9})^{pq} \le n^{p + 1.1q - 0.9pq}. $$
 So $\mathbb{E}[Z_{2,3}] \le n^{-0.1}$ and $\mathbb{E}[Z_{k,2}] \le n^{2.2 - 0.8k} \le n^{-0.2}$ for $k \ge 3$.
 % by Lemma~\ref{chernoff2},
Hence by Markov's inequality, $(ii)$ and $(iii)$ hold  with probability at least $1-o(1)$.

 Finally we show $(v)$.
 Suppose for all $x\in A$ and $v \in B$,
             $d_{H}(\{x,v\})\ge {n-1\choose
                   k-1}-{n-n/k\choose k-1}-\rho n^{k-1}$.
We see that, for $1 \le i \le n^{1.1}$ and for $x\in R_i \cap A$ and $v\in R_i \cap B$,
\begin{equation*}
\begin{split}
\mathbb{E}\left[d_{R_i}(\{x,v\})\right]
& > {n-1\choose k-1} n^{-0.9(k-1)}-{n-n/k\choose k-1}n^{-0.9(k-1)}-\rho n^{k-1}n^{-0.9(k-1)} \\
& > {n^{0.1}-1 \choose k-1} -{n^{0.1} - n^{0.1}/k \choose k-1}-\rho n^{0.1(k-1)}.
\end{split}
\end{equation*}
               %$$\mathbb{E}\left[d_{R_i}(\{x,v\})\right]> {n-1\choose k-1} n^{-0.9(k-1)}-{n-n/k\choose k-1}n^{-0.9(k-1)}-\rho n^{k-1}n^{-0.9(k-1)} > {n^{0.1} \choose k} -{n^{0.1} - mn^{-0.9} \choose k}-\rho n^{0.1k}.$$
By $(iv)$, with probability at least $1-O(n^{1.1})e^{-\Omega(n^{0.02})}$, for all $i= 1, \dots ,n ^{1.1}$, we have
              $|R_i\cap B| =(1+o(n^{-0.04}))n^{0.1}$.
Thus for all $x\in R_i \cap A$ and $v\in R_i \cap B$,
%there exists some $\rho'' >0$ such that
               $$\mathbb{E}\left[d_{R_i}(\{x,v\})\right] > {|R_i \cap B|\choose k-1}-{|R_i \cap B|-|R_i \cap B|/k \choose k-1}-2\rho |R_i \cap B|^{k-1}.$$
We wish to apply Lemma~\ref{janson} with $\Gamma = B$, $\Gamma_p = R_i \cap B$ and
$S={N_H(\{x,v\}) \cap B \choose k-1}$.
%$S$ be a family of all $(k-1)$-set of $B$ in $N_H(\{x,v\})$.
We define
$$\Delta = \frac{1}{2} \sum_{b_1, b_2 \in S, b_1 \ne b_2, b_1 \cap b_2 \ne \emptyset} \mathbb{E}[I_{b_1}I_{b_2}] \le \frac{1}{2} |R_i \cap B|^{2k-3} $$
%By Lemma~\ref{janson},
Thus,
\begin{align*}
&\mathbb{P}\left( d_{R_i}(\{x,v\}) \le {|R_i \cap B|\choose k-1}-{|R_i \cap B|-|R_i \cap B|/k \choose k-1}-3\rho |R_i \cap B|^{k-1} \right) \\
\leq & \mathbb{P} \left( d_{R_i}(\{x,v\}) \leq \mathbb{E}[d_{R_i}(\{x,v\})] - \rho |R_i \cap B|^{k-1} \right) \\
\leq & \exp(- \frac{(\rho |R_i \cap B|^{k-1})^2}{2{\mathbb{E}(d_{R_i}(\{x,v\}))} + 4\Delta}) \quad \mbox{(by Lemma~\ref{janson})}\\
\leq & \exp(- \frac{\rho^2 |R_i \cap B|^{2k-2}}{2{|R_i \cap B|\choose k-1} + 2|R_i \cap B|^{2k-3}}) \\
\leq & \exp(-\Omega(n^{0.1})).
\end{align*}
Therefore, with probability at least $1-O(n^{1.1})e^{-\Omega(n^{0.1})}$, $(v)$ holds.

By applying  union bound, $(i)$ -- $(v)$ all hold with probability $1-o(1)$.
\qed

\medskip

Now we prove that when the hypergraph $H$ in Theorem~\ref{general}
is not close to ${\cal H}_4(n,n/4)$ then $H$ contains an almost regular
spanning subgraph.

\begin{lemma}\label{Span-subgraph}
Let $0 < \rho \ll \varepsilon \ll 1$ be reals, and $n \equiv 0 \pmod 4$ be sufficiently large.
Let $H$ be a balanced $(1,4)$-partite $5$-graph with partition classes $X,[n]$ such that $|X|=n/4$.
Suppose that $d_{H}(\{x,v\})>{n-1\choose 3}-{3n/4 \choose 3}-\rho n^3$ for all $x\in X$ and $v \in [n]$.
If $H$ is not $\varepsilon$-close to $\mathcal{H}_4(n,n/4)$, then there exists a spanning subgraph $H'$ of $H$
such that the following conditions hold:
\begin{itemize}
		\item[$(1)$] For all $x\in V(H')$, with at most $n^{0.99}$ exceptions,
                               $d_{H'}(x)=(1\pm n^{-0.01})n^{0.2}$;
		\item[$(2)$] For all $x\in V(H')$, $d_{H'}(x)< 2 n^{0.2}$;
		\item[$(3)$] For any two distinct $x,y\in V(H')$, $d_{H'}(\{x,y\})< n^{0.19}$.
	\end{itemize}
\end{lemma}

\pf
Let $A_3\subseteq X$ and $A_4\subseteq [n]$ with $|A_i|=n^{0.99}$ for $i=3,4$.
Let $R_1,\ldots,R_{n^{1.1}}$ be defined as in Lemma \ref{lem1-5}. By $(iv)$ of Lemma \ref{lem1-5} , we have, for all $i= 1, \dots ,n ^{1.1}$,
              $$|R_i\cap X| =(1/4+o(n^{-0.04}))n^{0.1} \mbox{ and } |R_i\cap [n]| =(1+o(n^{-0.04}))n^{0.1}.$$
By $(v)$ of Lemma \ref{lem1-5}, we have,  for $1 \le i \le n^{1.1}$ and for $x\in X \cap R_i$ and $v \in [n] \cap R_i$,
               $$d_{R_i}(\{x,v\})>  {|R_i \cap [n]|\choose 3}-{3|R_i \cap [n]|/4 \choose 3}-3\rho |R_i \cap [n]|^{3};$$

By $(iv)$ and $(vi)$ of Lemma \ref{lem1-5},
we may choose $I_i\subseteq R_i\cap (A_3\cup A_4)$ such that for $i=1, \ldots, n^{1.1}$, $R_i' := R_i\setminus I_i$ is balanced and  $|R_i'|=(1 - o(1))|R_i|$.
%, where $R_i'=R_i\setminus I_i$ for $i=1, \ldots, n^{1.1}.$

%Since $H$ is not $\varepsilon$-close to $\mathcal{H}_t(k,n)$, $H_1$ is not $\varepsilon$-close to $\mathcal{F}_t(k,n)$ by Observation 2 in Section 2.
Let $a_1 = \varepsilon/32,  a_2=\varepsilon/512$, and $a_3 < \varepsilon(2^4 \cdot 4! \cdot 32)^{-1}$.
By  applying Lemma \ref{dense} to $H, a_1, a_2, a_3$,  we see that $H$ is $(\mathcal{F},a_1)$-dense, where
$$\mathcal{F}=\{U\subseteq V(H) : |U\cap X|\ge (1/4-a_1 ) n,\ |U\cap [n]|\ge (3/4-a_2 ) n\}.$$

%{\color{red}Delete ``Note that $e(H_1) \ge \beta e(H)$ where $\beta = (1 - (1-1/3k^2)^k)/3k \ge 1/9k^3$.
%Hence, $H$ is $(\mathcal{F},\varepsilon/54k^4)$-dense."}

Now we apply Lemma \ref{indep} to $H_1$ with $l = (3 \cdot 4^3 \cdot 4!)^{-1}$, $\alpha_1 = 1/4 - a_1$, $\alpha_2 = 3/4 -a_2$, and $\varepsilon' =a_1$.
Therefore, with probability at least $1-n^{O(1)}e^{-\Omega (n^{0.1})}$, for any independent set $S$ of $R_i'$,  $|S\cap R_i'\cap X | \leq (1/4 - a_1+o(1))n^{0.1}$ or $|S\cap R_i'\cap [n]|\leq (3/4-a_2+o(1))n^{0.1}$.

By applying Lemma \ref{frac-PM} to each $H[R_i']$, we see that each $H[R_i']$ contains a fractional perfect matching $\omega_i$.
Let $H^*=\cup_{i=1}^{n^{1.1}}R_i'$. We select a generalized binomial subgraph $H'$ of $H^*$ by letting $V(H')=V(H)$ and
independently choosing edge $e$ from $E(H^*)$,
with probability $\omega_{i_e}(e)$ if $e\subseteq R_{i_e}'$. (By  $(iii)$ of Lemma \ref{lem1-5}, for each $e\in E(H^*)$,  $i_e$ is uniquely defined.)

Note that since $w_i$ is a fractional perfect matching of $H[R_i']$ for $1 \le i \le n^{1.1}$, $\sum_{e \ni v} w_i(e) \le 1$ for $v \in R_i'.$
By $(i)$ of Lemma \ref{lem1-5}  and by Lemma~\ref{chernoff1},   $d_{H'}(v)=(1 \pm n^{-0.01}) n^{0.2}$ for any vertex $v\in V(H)\backslash (\cup_{i=1}^{n^{1.1}} I_i) \subseteq V(H)\backslash (A_3 \cup A_4)$
and $d_{H'}(v) \le (1 + n^{-0.01}) n^{0.2} < 2n^{0.2}$ for vertex $v \in \cup_{i=1}^{n^{1.1}} I_i$.
By $(ii)$ of Lemma \ref{lem1-5},  $d_{H'}(\{x,y\})\leq 2 < n^{0.19}$ for any $\{x,y\}\in {V(H) \choose 2}$.
%We choose $D=n^{0.2}$, $\tau=\tau'=n^{-0.01}$. Then
Therefore, $H'$ is the desired hypergraph.
\qed

\medskip

%% We may need to cite Pippenger instead of Frankl and Rodl.
We also need the following lemma attributed to Pippenger and Spencer \cite{PS} (see Theorem 4.7.1 in \cite{AS08}), which extends a result of Frankl and R\"odl \cite{FR85}.

\begin{lemma}[Pippenger and Spencer \cite{PS}, 1989] \label{nibble}
	For every integer $k\ge 2$ and reals $r> 1$ and $a>0$, there are $\gamma=\gamma(k,r,a)>0$ and $d_0=d_0(k,r,a)$ such that for every positive integer $n$ and $D\ge d_0$ the following holds: Every $k$-uniform hypergraph $H=(V,E)$ on a set $V$ of $n$ vertices in which all vertices have positive degrees and which satisfies the following conditions:
	\begin{itemize}
	 	\setlength{\itemsep}{0pt}
		\setlength{\parsep}{0pt}
		\setlength{\parskip}{0pt}
		\item[$($1$)$] For all vertices $x\in V$ but at most $\gamma n$ of them, $d(x)=(1\pm \gamma)D$;
		\item[$($2$)$] For all $x\in V$, $d(x)<rD$;
		\item[$($3$)$] For any two distinct $x,y\in V$, $d(x,y)<\gamma D$;
	\end{itemize}
	 contains a matching of size at least $(1-(k-1)a)(n/k)$.
\end{lemma}

\begin{comment}

\begin{lemma}\label{al-PM}
%For any $\varepsilon>0$ there exists $\rho, \rho'$ such that
$0<\rho'\ll \rho\ll \varepsilon\ll 1$.
Let $H$ be a  $(1,4)$-partite $5$-graph on the same vertex set with $\mathcal{H}_4(n,n/4)$. If $H$ is not $\varepsilon$-close to $\mathcal{H}_4(n,n/4)$,  and $d_H(\{u,v\})\geq {n-1\choose 3}-{3n/4\choose 3}-\rho n^3$ for any $v\in Q$ and $u\in [n]$, then $H$ has a matching covering all but at most $\rho'n$ vertices.
\end{lemma}

\pf
By Lemma \ref{F-R-subgraph}, $H$ has a spanning subgraph $H'$ such that $(1-\tau)D <deg_{H'}(v)<(1+\tau)D$ for all $v\in V$, and $\Delta_2(H') <\tau D$.
By Lemma \ref{regular}, $H'$ has a matching covering all but at most $\rho'n$ vertices. Thus $H$ has the desired matching. \qed

\end{comment}

\section{ \textbf{Proof of Theorem \ref{general} }}

\noindent {\bf Proof of Theorem~\ref{general}}.
Let $0<\rho'\ll \rho\ll \eta \ll \varepsilon\ll 1$. %Write $H:=H_{n,n/4}^4(\mathcal{F})$.
By Lemma \ref{close}, we may assume that $H$ is not $\varepsilon$-close to $\mathcal{H}_4(n,n/4)$. By Lemma \ref{absorbing}, there exists a matching $M_1$ with $|M_1|\leq \rho n$ such that for any balanced set $S$ with $|S|\leq \rho' n$, $H[V(M_1)\cup S]$ has a perfect matching.

Let $H_1=H-V(M_1)$. Then $H_1$ is not $(\varepsilon/2)$-close to $\mathcal{H}_4(n-4|M_1|,n/4-|M_1|)$. Write $n_1=|[n] \backslash V(M_1)|$. Furthermore, for all $x \in V(H_1) \cap X$ and $v \in V(H_1) \cap [n]$,
\[
d_{H_1}(\{x,v\})\geq {n-1\choose 3}-{3n/4\choose 3}-4|M_1|n^2>{n_1-1\choose 3}-{3n_1/4\choose 3}-10\rho n_1^3.
\]
%Let $n_1:=|H_1\cap [n]|$.
By Lemma \ref{Span-subgraph}, $H_1$ has a spanning subgraph $H_1'$  such that the following conditions hold:
\begin{itemize}
		\item[$(1)$] For all $x\in V(H_1')$, with at most $n_1^{0.99}$ exceptions,
                               $d_{H_1'}(x)=(1\pm n_1^{-0.01})n_1^{0.2}$;
		\item[$(2)$] For all $x\in V(H_1')$, $d_{H_1'}(x)< 2 n_1^{0.2}$;
		\item[$(3)$] For any two distinct $x,y\in V(H_1')$, $d_{H_1'}(\{x,y\})< n_1^{0.19}$.
	\end{itemize}
%$H_1$ has a spanning subgraph $H_1'$ such that $(1-\tau)D <d_{H'}(v)<(1+\tau)D$ for all $v\in V(H_1')$, and $\Delta_2(H') <\tau D$.
By Lemma \ref{nibble}, $H_1'$ has a matching, say $M_2$, covering all but at most $\rho'n$ vertices.
Write $S=V(H_1) \backslash V(M_2)$. Recall $H[V(M_1)\cup S]$ has a perfect matching $M_1'$.
Now $M_1'\cup M_2$ gives a desired perfect matching.
%This completes the proof.
\qed

%By Lemma \ref{al-PM}, $H_1$ has a matching covering $M_2$ all but at
%most $\rho'n$ vertices.

 %has no independent set $S$ such that with $|S\cap (Q-V(M_1))|\geq n_1/4-\varepsilon n_1$ and $|S\cap ([n]-V(M_1))|\geq n_1- 4\varepsilon n_1$. Note that
%\[
%d_{H_1}(\{u,v\})\geq {n-1\choose 3}-{3n/4\choose 3}-4\rho n^3\geq {n_1-1\choose 3}-{3n_1/4\choose 3}-5\rho n_1^3.
%\]
%By $H_1$ has a matching $M_2$ covering all but at most $\rho'n$ vertices.

\end{document}